\documentclass[12pt]{amsart}

\usepackage{vmargin}
\setmargrb{1in}{1in}{1in}{1in}
\usepackage{amstext}
\usepackage{amssymb,mathrsfs}
\usepackage{latexsym}
\newcommand\RR{{\mathbb R}}

\newcommand\NN{{\mathbb N}}
\newcommand\ZZ{{\mathbb Z}}

\newcommand\HH{{\mathbb H}}

\def\d={\,:=\,}

\font\frakten=eufm10
\newfam\frakfam
\textfont\frakfam=\frakten

\newtheorem{thm}{Theorem}[section]
\newtheorem{lemma}[thm]{Lemma}

\newtheorem{prop}[thm]{Proposition}
\newtheorem{Defn}[thm]{Definition}
\newtheorem{Ex}[thm]{Example}
\newtheorem{Rem}[thm]{Remark}
\newtheorem{Exs}[thm]{Examples}
\newtheorem{Rems}[thm]{Remarks}
\newtheorem{Defrem}[thm]{Definition and Remark}
\newtheorem{Remnt}[thm]{}

\newenvironment{rem}
 {\begin{Rem} \begin{rm}} {\end{rm} \hfill $\Box$ \end{Rem}}

\newenvironment{prf} {{\bf Proof.}}{\hfill $\Box$}

\begin{document}

\author{H.~F{\"u}hr
}

\address{Institute of Biomathematics and Biometry \\
GSF Research Center for Environment and Health \\
Ingolst{\"a}dter Landstra{\ss}e~1 \\
D-85764 Neuherberg \\
Germany.}
\email{fuehr@gsf.de}

\title[Vector-valued Gabor frames of Hermite functions]{Simultaneous estimates for vector-valued Gabor frames of Hermite functions}


\date{\today}

\begin{abstract}
We derive frame estimates for vector-valued Gabor systems with
window functions belonging to Schwartz space. The main result
provides frame bound estimates for windows composed of Hermite
functions. The proof is based on a recently established sampling
theorem for the simply connected Heisenberg group, which is
translated to a family of frame estimates via a direct integral
decomposition.
\end{abstract}

\maketitle

\section{Introduction and main results}  \label{sect:intro}

The aim of this paper is to derive frame estimates for
vector-valued Gabor systems. We consider the space ${\rm
L}^2(\mathbb{R};\mathbb{C}^d)$ of vector-valued signals. Elements
of this space can also be understood as vectors $\mathbf{f} =
(f_1,\ldots,f_d)$, with $f_i \in {\rm
L}^2(\mathbb{R};\mathbb{C})$, which amounts to identifying ${\rm
L}^2(\mathbb{R};\mathbb{C}^d)$ with the $d$-fold direct sum of
${\rm L}^2(\mathbb{R})$. Gabor systems in this space are obtained
by picking a {\em window function} $\mathbf{f} \in {\rm
L}^2(\mathbb{R};\mathbb{C}^d)$, and applying translations and
modulations. Here the translation and modulation operators are
given as
\[
(T_y \mathbf{f})(x) = \mathbf{f}(x-y)~~,~~(M_\xi \mathbf{f})(x) =
e^{2\pi i \xi x} \mathbf{f}(x) ~~(y,\xi \in \mathbb{R})
\]
 For any $\gamma =
(\gamma_1,\gamma_2) \in \mathbb{R}^2$, we denote the associated
time-frequency shift of a function $\mathbf{f} \in {\rm
L}^2(\mathbb{R};\mathbb{C}^d)$ by
\[
 \mathbf{f}_\gamma = T_{\gamma_1} M_{\gamma_2} \mathbf{f}~~.
\]

Now, given a lattice $\Gamma \subset \mathbb{R}^2$ and $\mathbf{f}
\in {\rm L}^2(\mathbb{R};\mathbb{C}^d)$, the resulting Gabor
system $\mathcal{G}(\mathbf{f},\Gamma)$ is given by
\[
 \mathcal{G}(\mathbf{f},\Gamma) = ( \mathbf{f}_\gamma )_{\gamma
 \in \Gamma}~~.
\]
Recall that a family $(\eta_i)_{i \in I}$ of vectors in a Hilbert
space $\mathcal{H}$ is called a frame if it satisfies
\begin{equation}
A \| \varphi \|^2 \le \sum_{i \in I} |\langle \varphi, \eta_i
\rangle |^2 \le B \| \varphi \|^2~~,
\end{equation} for all $\varphi \in \mathcal{H}$, with constants $0 < A \le B$. These
constants are called {\em frame bounds}; they are generally not
assumed to be optimal. A Gabor system
$\mathcal{G}(\mathbf{f},\Gamma)$ that is also a frame is called a
{\em Gabor frame}.

For the case $d=1$, Gabor frames have been studied extensively;
see e.g. \cite{Gr}. For a treatment of the case $d >1$, in
somewhat different terminology, see e.g. \cite{Ba}.

The problem of constructing a Gabor frame in dimension $d$
contains that of constructing $d$ Gabor frames: If
$\mathcal{G}(\mathbf{f},\Gamma)$ is a frame with bounds $A,B$,
then $\mathcal{G}(f_i,\Gamma)$ is a frame with frame bounds $A,B$,
for each component $f_i$ of $\mathbf{f}$. More generally,
$\mathcal{G}(\tilde{\mathbf{f}},\Gamma)$ is a frame of ${\rm
L}^2(\mathbb{R};\mathbb{C}^{\tilde{d}})$, for $\tilde{\mathbf{f}}
= (f_1,\ldots,f_{\tilde{d}})$, with $\tilde{d} < d$, again with
frame bounds $A$ and $B$.

But the converse need not be true: Note that by definition of the
scalar product in ${\rm L}^2(\mathbb{R};\mathbb{C}^d)$,
\[
 \langle \mathbf{g}, \mathbf{f}_\lambda \rangle = \sum_{i=1}^d
 \langle g_i, f_i|_\lambda \rangle ~~,
\] where we used $f_i|_\lambda$ to denote the action of a time-frequency shift
on $f_i$. Hence, cancellation may prevent the higher-dimensional
system from being a frame even when all components
$f_1,\ldots,f_d$ generate a frame. In fact, one can easily see
that a obvious necessary requirement for $\mathbf{f}$ to generate
a frame is linear independence of its entries $f_1,\ldots,f_d$.

Observe however, that at least the upper frame bound for
$\mathbf{f}$ can be estimated in terms of the upper frame bounds
for the $f_i$: If $B,B_1,\ldots,B_d$ are optimal upper framebounds
for $\mathbf{f},f_1,\ldots,f_d$, respectively, then the
Cauchy-Schwartz inequality entails \begin{equation}
\label{eqn:sf_fr_bd} B \le d \sum_{i=1}^d B_i ~~.\end{equation}

Probably the most-studied window function for the case $d=1$ has
been the Gaussian, $ g(x)= \pi^{-1/4} e^{- x^2/2}$. This is partly
due to historical reasons: Gabor suggested using Gaussian windows
\cite{Gabor}, and the characterization of densities for Gabor
frames with Gaussian window took more than 30 years to be fully
clarified \cite{Lyub,SeipWallsten}. The choice of this window
function is motivated by the way Gabor systems are employed: They
are designed to measure time-frequency content in a signal. By the
Heisenberg uncertainty relation a Gaussian window has optimal
time-frequency concentration, and thus can be expected to yield a
good time-frequency resolution. Moreover, for the Gaussian window
powerful tools from complex analysis can be employed to study
sampling \cite{Lyub,SeipWallsten}, which adds to its theoretical
appeal.

In this paper, we intend to derive frame estimates for window
functions $\mathbf{f}$ consisting of the first $d+1$ Hermite
functions. For $n \in \NN_0$, we define the $n^{\rm th}$ Hermite
function $h_n$ by
 \begin{equation}
  h_n(x) = \frac{(-1)^n}{\sqrt{2^n n ! \sqrt{\pi}}} ~ e^{x^2/2} ~ \frac{d^n}{dx^n}e^{-x^2}~~,
 \end{equation} where the normalization factor ensures $\| h_n \|_2 = 1$.
The above defined Gaussian equals $h_0$, whence the problem
considered here can be viewed as a generalization of Gabor's
original question.

The vector-valued windows that we are interested in are given by
\[
 \mathbf{h}^d = \left( h_0,\ldots, h_d \right) \in {\rm
 L}^2(\mathbb{R}; \mathbb{C}^{d+1})~~.
\] We intend to give frame conditions and estimates for
$\mathcal{G}(\mathbf{h}^d,\mathcal{M}(\mathbb{Z}^2))$, with
$\mathcal{M}$ an invertible matrix, in terms of a matrix norm
defined by \begin{equation} \label{eqn:defn_matr_norm} \|
\mathcal{M} \| = \sup \{ \|\mathcal{M} z \|_2 : \|z \|_\infty \le
1/2 \} ~~.\end{equation} Here $\|z\|_p$ denotes the usual
$\ell^p$-norm on $\RR^2$. This choice of matrix norm may seem
somewhat peculiar, and in fact the theorems below can be
formulated with respect to any other norm on matrix space. The use
of (\ref{eqn:defn_matr_norm}) emphasizes the close connection to
sampling estimates on the Heisenberg group.

Let us now state the main results of this paper. The chief purpose
of Theorem \ref{thm:main1} is to allow a better understanding and
formulation of Theorem \ref{thm:main2}. Nonetheless, Theorem
\ref{thm:main1} is of independent interest. Even though we expect
it somehow to be part of Gabor analysis folklore, we are not aware
of any previous source for this result; not even for $d=1$.

\begin{thm} \label{thm:main1}
 Let $\mathbf{f} = (f_1,\ldots,f_d) \in {\rm L}^2(\mathbb{R}, \mathbb{C}^d)$ be given,
 with $f_i \in \mathcal{S}(\mathbb{R})$, and $\langle f_i,f_j \rangle = \delta_{i,j}$, for all
 $1 \le i,j \le d$. Then there exists a constant $0 < C_{\mathbf{f}} \le 1$ such that for all
 matrices $\mathcal{M}$ with $\|\mathcal{M} \| < C_{\mathbf{f}}$, the system
 $\mathcal{G}(\mathbf{f},\mathcal{M}(\mathbb{Z}^2))$ is a frame with frame constants $\frac{1}{|\det(A)|} \left( 1 \mp \frac{\|
 {\mathcal{M}} \|}{C_{\mathbf{f}}}
 \right)^2$.
\end{thm}

Hence the {\em tightness} of the frame estimate, which is the
quotient of the two frame bounds, approaches $1$ as $\|\mathcal{M}
\| \to 0$, with speed proportional to $\| \mathcal{M} \|^2
/C_{\mathbf{f}}^2$.

We observe that Theorem \ref{thm:main1} holds for the supremum
$C_{\mathbf{f}}^*$ of all possible constants, and that this choice
provides the sharpest possible statement. Then the main result of
this paper is the following estimate:

\begin{thm} \label{thm:main2} There exists a constant $C_{\mathbb{H}} \le 1$ such that for all $d \in \mathbb{N}_0$
\begin{equation} C_{\mathbf{h}^d}^* \ge
\frac{C_{\mathbb{H}}}{\sqrt{2d+1}}~~. \end{equation}
\end{thm}

\section{Method of proof: Sampling vs. frame estimates}

The proof of Theorem \ref{thm:main2} relies on a combination of
various techniques: It uses a sampling estimate for the
Paley-Wiener space $PW(\mathbb{H})$, established in \cite{FuGr}.
The definition of this space uses a particular differential
operator on $\mathbb{H}$, the so-called sub-Laplacian. Hermite
functions enter in the spectral decomposition of this operator,
and it is this connection that will allow to relate the sampling
theorem to frame estimates for Hermite functions.

The connection between frames and sampling theory  is not exactly
new, in fact it is at the base of frame theory, which originated
from nonharmonic Fourier series and their connections to irregular
sampling over the reals, see the classic paper \cite{DS52}. For
the sake of explicitness, assume we are given a sequence $\Lambda
= (\lambda_k)_{k \in \ZZ}$ of sampling points in $\RR$. We are
looking for criteria that allow to reconstruct Paley-Wiener
functions, i.e. $f \in {\rm L}^2(\RR)$ whose Fourier transform has
support in the unit interval $[-0.5,0.5]$, in a stable manner from
its restriction to $\Lambda$. Noting that
\[ f (\lambda_k) = \int_{-0.5}^{0.5} \widehat{f}(\xi) e^{2 \pi i
\lambda_k \xi} d\xi = \langle \widehat{f}, e_{\lambda_k} \rangle
~~,\] we find the following two equivalent conditions, with
identical constants $A$ and $B$ in both cases:
\begin{enumerate}
\item The sequence $\Lambda$ fulfills the {\em sampling estimate}
\[ A \| f \|_2^2 \le \sum_{k \in \ZZ} |f(\lambda_k)|^2 \le B \| f
\|_2^2 ~~,\] for all Paley-Wiener functions. \item The sequence
$(e_{\lambda_k})_{k \in \ZZ}$ fulfills the frame estimate
\[
 A \| F \|_2^2 \le \sum_{k \in \ZZ} |\langle F, e_{\lambda_k}
 \rangle |^2 \le B \| F \|_2^2 ~~,
\] for all $F \in {\rm L}^2([-0.5,0.5])$.
\end{enumerate}

This equivalence can be used in two ways: For instance, observing
that the choice $\lambda_k = k$ ($k \in \ZZ$) results in the
Fourier orthonormal basis $(e_k)_{k \in \ZZ}$ of ${\rm
L}^2([-0.5,0.5])$, the implication $(2) \Rightarrow (1)$ leads to
Shannon's sampling theorem. Conversely, for irregular sampling
sets, condition (1) can often be checked using tools from complex
analysis techniques, and then $(1) \Rightarrow (2)$ results in
frame estimates for irregularly spaced exponentials.

In this paper, we use a similar approach for the Heisenberg group
$\HH$: This time, previously established sampling estimates for
Paley-Wiener functions on $\HH$ will allow to derive frame
estimates for Hermite functions, by an analogue of the implication
$(1) \Rightarrow (2)$. For this purpose we will need to work out
the connections between $PW({\mathbb{H}})$ and the Hermite
functions, which is the topic of Sections 4 and 5. But first let
us prove Theorem \ref{thm:main1}.

\section{Proof of Theorem \ref{thm:main1}}

 Let $V_{\mathbf{f}} : {\rm L}^2(\mathbb{R};\mathbb{C}^d) \to {\rm L}^2(\mathbb{R}^2)$ denote the
 windowed Fourier transform,
 \[
  V_{\mathbf{f}} \mathbf{g} (x,\xi) = \langle \mathbf{g}, T_x M_\xi \mathbf{f} \rangle~~.
 \]
 The orthogonality relations for the windowed Fourier transform and the pairwise orthogonality
 of the components of $\mathbf{f}$ imply that
 $V_{\mathbf{f}}$ is an isometry. Hence its image $\mathcal{H}_{\mathbf{f}}$ is a closed subspace of ${\rm
 L}^2(\mathbb{R}^2)$. As outlined in the previous subsection, the
 isometry property of $V_{\mathbf{f}}$ implies that a frame estimate
 for $\mathcal{G}(\mathbf{f},\mathcal{M}(\mathbb{Z}^2))$ is the same as a sampling
 estimate for $\mathcal{H}_{\mathbf{f}}$, with sampling set $\mathcal{M}(\mathbb{Z}^2)$.

 We intend to utilize the techniques from \cite{FuGr} for this purpose, hence we will need oscillation
 estimates. We define
 \[ {\rm osc}_r(f)(x) = \sup_{|x-y|<r} |f(x)-f(y)| ~~.\]
 Our first aim is to prove
 \begin{equation} \label{eqn:osc_Hf} \| {\rm osc}_r(F) \|_2 \le \frac{r}{C_{\mathbf{f}}} \| F \|_2  ~~,
 \forall F \in \mathcal{H}_{\mathbf{f}} ~~,\end{equation} for a
 suitable constant $C_{\mathbf{f}}$.

 For this purpose we first observe that the projection
 $P_{\mathbf{f}}$ onto $\mathcal{H}_{\mathbf{f}}$ is obtained by twisted
 convolution \cite{Gr}:
 \begin{eqnarray*}
  (P_{\mathbf{f}} G)(x,\xi) & = & (G \sharp F)(x,\xi) \\ & = &
  \int_{\mathbb{R}^2} G(x',\xi') F (x-x',\xi-\xi') e^{\pi i (x \xi'-x' \xi)}
  dx' d\xi'~~,
 \end{eqnarray*}
where we let $F = V_{\mathbf{f}} \mathbf{f}$.

Hence $G = G \sharp F$ for $G \in {\mathcal{H}}_{\mathbf{f}}$, and
therefore
\begin{eqnarray*}
\lefteqn{ osc_r(G)(x,\xi) } \\ & = & \sup_{|(x,\xi) -
(x'',\xi'')|<r} \left| \int_{\mathbb{R}^2}  G(x',\xi')
  \left( F (x-x',\xi-\xi') -  F (x''-x',\xi''-\xi') \right) e^{\pi i (x \xi'-x'
  \xi)}
  dx' d\xi' \right| \\
  & \le & \int_{\mathbb{R}^2}  |G(x',\xi')|  \sup_{|(x,\xi) - (x'',\xi'')|<r}
  \left| F (x-x',\xi-\xi') -  F(x''-x',\xi''-\xi') \right|
  dx' d\xi' \\
& = & \int_{\mathbb{R}^2}  |G(x',\xi')|  {\rm osc}_r(F)(x-x',\xi-\xi')  dx' d\xi' \\
& = & |G| \ast {\rm osc}_r (F) ~~,
\end{eqnarray*}
where convolution is taken with reference to the abelian group
structure on $\mathbb{R}^2$. Hence
\begin{equation}
 \label{eqn:osc_twist} \| {\rm osc}_r (G) \|_2 \le \| G \|_2 \| {\rm osc}_r (F) \|_1~~.
\end{equation}

Now the second factor can be estimated by
\begin{eqnarray*}
\| {\rm osc}_r (F) \|_1 & = & \int_{\mathbb{R}^2} \sup_{(x',\xi') \in B_r(x,\xi)} |F(x,\xi)-F(x',\xi')| dxd\xi \\
& \le & \int_{\mathbb{R}^2} r \sup_{(x',\xi') \in B_r(x,\xi)}
\left( \left| \frac{\partial F}{\partial x}(x,\xi) \right| +
 \left| \frac{\partial F}{\partial \xi}(x,\xi) \right| \right) dx d\xi \\
& \le & r \int_{\mathbb{R}^2} \sum_{|\alpha| \le 1} \|D^\alpha F \|_{\infty,B_r(x,\xi)} dx d\xi \\
& \le & r \int_{\mathbb{R}^2} \sum_{|\alpha| \le 4} \|D^\alpha F \|_{1,B_r(x,\xi)} dx d\xi \\
& = & r \sum_{|\alpha|\le 4} |B_r| C_{B_r} \| D^\alpha F \|_1
\end{eqnarray*}
where we used the mean value theorem for the first inequality, and
the Sobolev embedding theorem for the third. Here $C_{B_r}$
denotes the norm of the embedding $W^{3,1}(B_r) \hookrightarrow
C(B_r)$ \cite[Theorem 5.4]{Adams}. Clearly, $|B_r| = r^2 |B_1|$;
on the other hand, a dilation argument establishes for $r<1$ that
$C_{B_r} \le r^{-2} C_{B_1}$. Hence
\[ \| {\rm osc}_r (F) \|_1 \le r |B_1| C_{B_1} \sum_{|\alpha|\le 4} \| D^\alpha F
\|_1~~,\] and thus (\ref{eqn:osc_twist}) implies
(\ref{eqn:osc_Hf}), with $C_f = \frac{\displaystyle
1}{\displaystyle |B_1| C_{B_1} \sum_{|\alpha|\le 4} \| D^\alpha F
\|_1}$.

 Now by letting $K = \mathcal{M}([-0.5,0.5]^2)$, we obtain $\RR^2 =
 \bigcup_{\gamma \in \mathcal{M}(\mathbb{Z}^2)} \gamma + K$ as a
 disjoint union, and $K$ has Lebesgue measure $| {\rm det}(\mathcal{M})|$.
 Moreover, by definition of the norm, $K \subset
 B_r$, for $r = \| \mathcal{M} \|$. Now, by \cite[Theorem
 3.5]{FuGr}, the oscillation estimate (\ref{eqn:osc_Hf})
 results in the sampling estimate
 \[
  \frac{1}{| {\rm det}(\mathcal{M})|} \left( 1- \frac{r}{C_{\mathbf{f}}}
  \right)^2 \| F \|_2^2 \le \sum_{\gamma \in
  \mathcal{M}(\mathbb{Z}^2)} |F(\gamma)|^2 \le  \frac{1}{| {\rm det}(\mathcal{M})|}
  \left( 1+ \frac{r}{C_{\mathbf{f}}} \right)^2 \| F \|_2^2~~,
 \] which is Theorem \ref{thm:main1}.
\hfill $\Box$

\begin{rem}
 As the proof shows, one can weaken the requirement
 $f_i \in \mathcal{S}(\mathbb{R})$ to the membership of
 $V_{\mathbf{f}} \mathbf{f}$ in a suitable Sobolev space.
 In addition, pairwise orthogonality of the $f_i$ can be weakened to
 linear independence; in this case the frame estimate
 will contain a term involving the Gramian matrix of the $f_i$.
\end{rem}

\section{Fourier transform on the Heisenberg group}

The (simply connected) Heisenberg group is defined as $\mathbb{H}
= \mathbb{R}^3$, with group law
\[ (p,q,t) (p',q',t') = (p+p',q+q',t+t'+(pq'-p'q)/2) ~~.\] For the following facts concerning
$\mathbb{H}$, we refer the reader to \cite{Fo}. $\mathbb{H}$ is a
step-two nilpotent Lie group, with center $Z(\mathbb{H}) = \{ 0 \}
\times \{ 0 \} \times \RR$. $\mathbb{H}$ is unimodular, with
two-sided invariant measure on $\mathbb{H}$ given by the usual
Lebesgue measure of $\mathbb{R}^3$.

Given $\lambda \in \RR^* = \RR \setminus \{ 0 \}$, the {\em
Schr\"odinger representation} $\rho_\lambda$ of $\mathbb{H}$ acts
on ${\rm L}^2(\mathbb{R})$ via
\[
\rho_\lambda (p,q,t) = e^{-2 \pi i \lambda (t-pq/2)} T_{\lambda p}
M_q ~~.
\]
This is an irreducible unitary representation of $\mathbb{H}$. The
family of Schr\"odinger representations provides the basis for the
Plancherel transform of the group, a tool that is of key
importance for this paper.

 Before we describe this transform in more detail,
let us quickly recall the basics of Hilbert-Schmidt operators: The
space of Hilbert-Schmidt operators on a Hilbert space
${\mathcal{H}}$ is given by all bounded linear operators such that
\[ \| T \|_{HS}^2 = \sum_{i \in I} \| T \eta_i \|^2 \]
is finite; here $(\eta_i)_{i \in I}$ denotes an arbitrary
orthonormal basis of ${\mathcal{H}}$. It is well-known that the
norm is independent of the choice of basis, and defines a Hilbert
space structure on the Hilbert-Schmidt operators with scalar
product
\[ \langle S, T \rangle = {\rm trace}(T^*S) = \sum_{i \in I}
\langle S \eta_i, T \eta_i \rangle~~.
\]

Now, given $f \in {\rm L}^1 \cap L^2(\mathbb{H})$, one defines
\[ \rho_\lambda(f) = \int_{\mathbb{H}} \rho_\lambda(x) f(x) dx~~, \]
understood in the weak operator sense. This is just the canonical
extension of the representation $\rho_\lambda$ to the convolution
algebra ${\rm L}^1({\mathbb{H}})$. The mapping $f \mapsto
(\rho_\lambda(f))_{\lambda \in \mathbb{R}^*}$ is called the {\em
group Fourier transform}. This nomenclature is justified by the
observation that the euclidian Fourier transform is obtained by
integration against the characters of the additive group.
Moreover, it turns out that for $f \in {\rm L}^1 \cap
L^2(\mathbb{H})$, the group Fourier transform is in fact a family
of Hilbert-Schmidt operators, and we have the {\em Parseval
relation}
\[
\| f \|_2^2 = \int_{\RR} \| \rho_\lambda(f) \|_{HS}^2 |\lambda|
d\lambda~~.
\]

The Parseval relation allows to extend the Fourier transform to
${\rm L}^2(\mathbb{H})$, yielding the {\em Plancherel transform},
a unitary map
\[ {\rm L}^2(\mathbb{H}) \to \int_{\mathbb{R}}^\oplus HS({\rm
L}^2(\mathbb{R})) ~|\lambda|d\lambda ~~,\] where the right hand
side denotes the direct integral of Hilbert-Schmidt spaces. We
denote the Plancherel transform of $f \in {\rm L}^2(\mathbb{H})$
as $\widehat{f} = \left( \widehat{f}(\lambda) \right)_{\lambda \in
\mathbb{R}^*}$. This map will play the same role as the euclidean
Fourier transform in the discussion of Section 2.

Of crucial importance for this usage of the Plancherel transform
are its algebraic properties, providing a decomposition of various
operators and representations acting on ${\rm L}^2(\mathbb{H})$.
For instance, if we denote the {\em left regular representation}
by $L_x f (y) = f(x^{-1} y)$, for $f \in {\rm L}^2(\mathbb{H})$
and $x,y \in \mathbb{H}$, then
\[ \rho_\lambda(L_x f) = \rho_\lambda(x) \circ \rho_\lambda(f) ~~,\]
which provides the decomposition of the left regular
representation $L$ into a direct integral,
\[ L \simeq \int_{\mathbb{R}^*} \rho_\lambda \otimes 1 ~|\lambda|d\lambda ~~.\]
Similarly, the right regular representation $R$ decomposes by the
formula $\rho_\lambda(R_x f) = \rho_\lambda(f) \circ
\rho_\lambda(x)^*$.

The decompositions extend to commuting operators: For any bounded
operator $T$ commuting with $L$, there exists a measurable field
$(\widehat{T}_\lambda)_{\lambda \in \mathbb{R}^*}$ of bounded
operators on ${\rm L}^2(\RR)$ satisfying $\rho_\lambda(T f) =
\rho_\lambda(f) \circ \widehat{T}_\lambda$, or in direct integral
notation
\[ T \simeq \int_{\mathbb{R}^*} {\rm Id} \otimes
\widehat{T}_\lambda ~|\lambda|d\lambda ~~.\]

\section{Paley-Wiener space on the Heisenberg group}

In this section we outline the definition of Paley-Wiener space on
$\mathbb{H}$ and its relation to Hermite functions. The central
role of Hermite function in the decomposition of the sub-Laplacian
has been observed previously, e.g. in \cite{Ge}; the results
presented below can be found also in \cite{Tha}.

We define a left-invariant differential operator $P$ on
$\mathbb{H}$ by \begin{equation} \label{eqn:def_P} (Pf)(p,q,t) =
\lim_{h \to 0} \frac{f((p,q,t) (h,0,0))- f(p,q,t)}{h}
~~,\end{equation} corresponding to the subgroup $\mathbb{R} \times
\{ 0 \} \times \{ 0 \}$, and $Q$ is a left-invariant operator
associated to $\{ 0 \} \times \mathbb{R} \times \{ 0 \}$ in the
same manner. $P,Q$ are viewed as elements of the Lie algebra
$\mathfrak{h}$ of $\HH$; we have $[P,Q] = T$, the infinitesimal
generator of the group center. This observation exhibits
$\mathfrak{h}$ as a {\em stratified} Lie algebra,
\[ \mathfrak{h} = V_1 + V_2 \] with $V_1 = {\rm span} (P,Q)$, and
$V_2 = \RR \cdot T = [V_1,V_1]$.

Of particular interest for analysis on these groups is the {\em
sub-Laplacian}; as the name suggests, it can be viewed as a
replacement for the Laplacian over $\mathbb{R}^n$. For the
Heisenberg group, this operator is defined by
\[
 \mathcal{L} = -P^2 - \frac{Q^2}{4 \pi^2}~~.
\]
The normalization of $Q$ is chosen for the sake of convenience.
$\mathcal{L}$ is a left-invariant positive unbounded operator on
${\rm L}^2(\mathbb{H})$. We denote its spectral measure by
$\Pi_{\mathcal{L}}$. The Paley-Wiener space on $\mathbb{H}$ is
then defined as \[ PW({\mathbb{H}}) = \Pi_{\mathcal{L}}([0,1])
({\rm L}^2(\mathbb{H})) ~~.\] Note that, up to normalization, this
definition is completely analogous to the definition of
bandlimited functions on $\RR$, since the euclidean Fourier
transform can also be read as the spectral decomposition of the
Laplacian. The projection $\Pi_{\mathcal{L}}([0,1])$ is
left-invariant, and is therefore decomposed by the group Fourier
transform into a direct integral. The following lemma provides an
explicit calculation of this decomposition via Hermite functions.
In the following, we use the notation $D_a \mathbf{f} (x) =
|a|^{-1/2} \mathbf{f}(|a|^{-1} x)$. As with translation and
modulation operators, we use the same symbol for operators acting
on scalar- and on vector-valued functions.

\begin{lemma}
\begin{enumerate}
\item[(a)] $(h_n)_{n \in \mathbb{N}_0}$ is an orthonormal basis of
${\rm L}^2(\mathbb{R})$. \item[(b)] The system $(h_n)_{n \in
\mathbb{N}_0}$ is an eigenbasis of the Hermite operator $Hf(x) =
x^2f(x) - f''(x)$,
\begin{equation}
 H h_n = (2n+1) h_n~~.
\end{equation}
\item[(c)] For every real $a \not= 0$, the dilated system
$(h_{n,a})_{n \in \mathbb{N}_0}$, defined by \[ h_{n,a}(x) =
(D_{|a|^{1/2}} h_n)(x) = |a|^{-1/4} h_n(|a|^{-1/2} x)~~,\] is an
eigenbasis of the scaled Hermite operator $H_a f(x) = x^2f(x) -
a^2 f''(x)$,
\begin{equation} \label{eqn:eigenHermite}
 H_a h_{n,a} = |a| (2n+1) h_{n,a}
\end{equation}
\item[(c)] The sub-Laplacian decomposes into a direct integral of
scaled Hermite operators:
\[ \rho_\lambda(\mathcal{L} f) = \rho_\lambda(f) \circ H_\lambda ~~,\]
for all $f \in C_c^\infty(\mathbf{H})$.
\end{enumerate}
\end{lemma}
\begin{prf} For part (a) confer \cite[Corollary 6.2, Theorem 6.14]{Fo_Four}. Part
(b) follows from this by straightforward computation. Part (c) is
established by formal calculation from (\ref{eqn:def_P}) and the
analogous formula for $Q$, using the decomposition of the right
regular representation.
\end{prf}

Parts (b) and (c) contain the ingredients of the direct integral
decomposition of $PW(\mathbb{H})$. For the precise formulation of
this result and its proof, the tensor product notation for
Hilbert-Schmidt operators will be useful. Given vectors $\eta,
\varphi$ in a Hilbert space $\mathcal{H}$, we let
\[ \eta \otimes \varphi : z \mapsto \langle z, \varphi \rangle
\eta ~~,\] which is a rank-one operator on $\mathcal{H}$. Note
that the notation is only conjugate linear in $\varphi$. The
Hilbert-Schmidt scalar product of two elementary tensors is
\[ \langle \eta \otimes \varphi, \eta' \otimes \varphi'
\rangle_{HS} = \langle \eta, \eta' \rangle_{\mathcal{H}} \langle
\varphi', \varphi \rangle_{\mathcal{H}} ~~.\] Moreover, for any
pair $S,T$ of bounded operators, $S \circ (\eta \otimes \varphi)
\circ T = (S \eta) \otimes (T^* \varphi)$.

Now, given any orthonormal basis $(\varphi_i)_{i \in I}$ of
$\mathcal{H}$, every Hilbert-Schmidt operator $T$ has a unique
decomposition
\[
 T = \sum_{i \in I} \eta_i \otimes \varphi_i ~~.
\]
Hence, if $S = \sum_{i \in I} \psi_i \otimes \varphi_i$ is another
Hilbert-Schmidt operator, we obtain for the scalar product
\begin{equation} \label{eqn:scal_HS}
 \langle T, S \rangle = \sum_{i \in I} \langle \eta_i, \psi_i
 \rangle~~.
\end{equation}

Observe in the formulation of the following proposition that
$\widehat{P}_\lambda$ involves the first $d(\lambda)+1$ Hermite
functions. For $|\lambda|>1$, we have $d(\lambda) = -1$, and thus
$\widehat{P}_\lambda = 0$.

\begin{prop}
 Letting
 \[
 d(\lambda) = \left\lfloor \frac{1}{2 |\lambda|} - \frac{1}{2} \right\rfloor
\]
and
 \[
  \widehat{P}_\lambda = \sum_{n=0, \ldots, d(\lambda)} h_{n,\lambda}
  \otimes h_{n,\lambda}~~,
 \]
 the projection onto Paley-Wiener space is given by
 \begin{equation} \label{eqn:def_P_hut} \left( \Pi_{\mathcal{L}}([0,1]) (f) \right)^\wedge(\lambda) =
 \widehat{f}(\lambda) \circ \widehat{P}_\lambda~~,~~ \forall f \in {\rm L}^2(\mathbb{H}) ~~.
 \end{equation}
 The operator field $(\widehat{P}_\lambda)_{\lambda \in \mathbb{R}^*}$
 is the Plancherel transform of a function $p
 \in {\rm L}^2(\mathbb{H})$, whence $\Pi_{\mathcal{L}}([0,1])(f) = f \ast
 p$.
 \end{prop}

 \begin{prf}
We apply the above considerations to the case $\mathcal{H} = {\rm
L}^2(\RR)$ and its orthonormal basis $(h_{n,\lambda})_{n \in
\mathbb{N}_0}$. Hence each Hilbert-Schmidt operator $T$ on ${\rm
L}^2(\mathbb{R})$ has a decomposition
\begin{equation}
 T = \sum_{n \in \mathbb{N}_0} \eta_n \otimes
 h_{n,\lambda}~~,
\end{equation} and we obtain from (\ref{eqn:eigenHermite}) that
\[
 T \circ H_\lambda = \sum_{n \in \mathbb{N}_0} |\lambda| (2n+1) \eta_n \otimes
 h_{n,\lambda}~~.
\] This shows that the map $T \mapsto \eta_n \otimes
h_{n,\lambda}$ can be understood as a projection onto an
eigenspace of the operator $T \mapsto T \circ H_\lambda$, with
associated eigenvalue $|\lambda| (2n+1)$. By definition of
Paley-Wiener space, only eigenvalues $\le 1$ are admitted, which
shows that the definition of $\widehat{P}_\lambda$ indeed yields
(\ref{eqn:def_P_hut}).

For the second statement, we compute the norm of the operator
field in the direct integral space. First observe that
$\widehat{P}_\lambda = 0$ for $|\lambda|>1$. Moreover, the squared
Hilbert-Schmidt norm of a projection equals its rank, whence $\|
\widehat{P}_\lambda \|_{HS}^2 = d(\lambda) < \frac{1}{2
|\lambda|}$, and thus
\[
\int_{\mathbb{R}^*} \| \widehat{P}_\lambda \|_{HS}^2 ~|\lambda|
d\lambda < \int_{-1}^{1} \frac{1}{2 |\lambda|} ~|\lambda| d\lambda
= 1 ~~.\] Hence $(\widehat{P}_\lambda)_\lambda$ has a preimage $p$
under the Plancherel transform. Finally, (\ref{eqn:def_P_hut}) and
the convolution theorem \cite[Theorem 4.18]{Fu_LN} provide that
$Pf = f \ast p$.
\end{prf}

The motivation for considering $PW(\mathbb{H})$ is the existence
of sampling estimates. The formulation of the sampling theorem
requires some additional notation. We fix a quasi-norm $| \cdot |
: {\mathbb{H}} \to \mathbb{R}^+_0$ by
\[
 |(p,q,t)| = (p^2 + q^2 + |t|)^{1/2}~~,
\] and write $B_r$ for the unit ball around $0$.
A discrete subset $\Gamma \subset {\mathbb{H}}$ is called a {\em
quasi-lattice} if there exists a relatively compact set $K \subset
{\mathbb{H}}$ such that $\mathbb{H} = \bigcup_{\gamma \in \Gamma}
\gamma K$, as a disjoint union. Such a set $K$ is called {\em
complement} of $\Gamma$.
\begin{thm} \cite[Theorem 5.11]{FuGr} \label{thm:sampl}
 There exists a constant $0<C_{\mathbb{H}} \le 1$ with the following property: For all quasi-lattices
 $\Gamma$ possessing a complement $K$ contained in a ball of radius $r < C_{\mathbb{H}}$ and
 all $f \in PW({\mathbb{H}})$
 \begin{equation} \label{eqn:samp_est}
 \frac{1}{|K|} (1-r/C_{\mathbb{H}})^2 \| f \|_2^2 \le \sum_{\gamma \in \Gamma} | f(\gamma) |^2
 \le \frac{1}{|K|} (1+r/C_{\mathbb{H}})^2 \| f \|_2^2~~.
 \end{equation}
\end{thm}

\section{Proof of Theorem \ref{thm:main2}}

We will now derive Theorem \ref{thm:main2} from Theorem
\ref{thm:sampl}, basically by explicit calculation. The following
lemma can be seen as an analog of $(1) \Leftrightarrow (2)$ from
Section 2. A version of this result was obtained in
\cite[Proposition 6.11]{Fu_LN}.

\begin{lemma} \label{lem:WH_transfer}
Suppose that $\Gamma' \subset \mathbb{H}$ is of the form $\Gamma'
= \Gamma \times \alpha \ZZ$, with $\Gamma \subset \RR^2$ and
$\alpha
>0$.  Consider the following statements:
\begin{enumerate}
\item[(a)] For all $f \in PW({\mathbb{H}})$,
\begin{equation} \label{eqn:samp_trans}
 A \| f \|_2^2 \le \sum_{\gamma \in \Gamma'} | f(\gamma) |^2 \le B  \| f \|_2^2~~.
\end{equation}
\item[(b)] For all $d \in \mathbb{N}_0$, and for almost all
$\lambda$ with $|\lambda| < \frac{1}{2d+1}$, the system
$\mathcal{G}(\mathbf{h}^d,|\lambda|^{1/2} \Gamma)$ is a frame of
${\rm L}^2(\mathbb{R}; \mathbb{C}^{d+1})$ with frame bounds $
\alpha |\lambda|^{-1} A$ and $\alpha |\lambda|^{-1} B$.
\end{enumerate}
Then $(a) \Rightarrow (b)$, and if $\alpha<1/2$, $(b) \Rightarrow
(a)$.

Moreover, if $\Gamma = \mathcal{M} (\mathbb{Z}^2)$, for a suitable
invertible matrix, the frame estimates in $(b)$ are valid for all
$\lambda < \frac{1}{2d+1}$.
\end{lemma}

\begin{prf}
For the proof of $(a) \Rightarrow (b)$, let $f \in PW(\mathbb{H})$
be given. Then we can write
\begin{equation} \widehat{f}(\lambda) = \sum_{i=0}^{d(\lambda)}
 \varphi_{i,\lambda} \otimes h_{i, \lambda}~~,
\end{equation} for suitable functions $\varphi_{i,\lambda} \in {\rm
L}^2(\mathbb{R})$. In the following, we also use the notations
\[ \Phi_\lambda =
\left( \varphi_{0,\lambda}, \ldots, \varphi_{d(\lambda), \lambda}
\right) \] and
\[ \mathbf{h}_\lambda = \left( h_{0,\lambda}, \ldots,
h_{d(\lambda),\lambda} \right) ~~.\]

Let $E \subset{\mathbb{R}^*}$ be a Borel set contained in an
interval $I$ of length $1/\alpha$, and consider $g$ with
$\widehat{g} = \widehat{f} \cdot \mathbf{1}_E$. Observing that $g
\in PW(\mathbb{H})$, we can compute the $\ell^2$-norm of its
restriction to $\Gamma'$ as follows:
\begin{eqnarray*}
\sum_{\gamma \in \Gamma'} | g(\gamma) |^2 & = & \sum_{\gamma \in
\Gamma'} | \langle g, L_{\gamma} p \rangle |^2 \\ & = &
\sum_{\gamma \in \Gamma'} | \langle \widehat{g}, \left( L_{\gamma}
p\right)^\wedge \rangle |^2 \\ & = & \sum_{\gamma \in \Gamma'}
\left| \int_{E} \langle \widehat{f} (\lambda), \rho_\lambda
(\gamma)
\widehat{p}(\gamma) \rangle |\lambda| d\lambda \right|^2 \\
& = & \sum_{(l,k) \in \Gamma,n \in \ZZ}  \left| \int_E \langle \widehat{f} (\lambda),
\rho_\lambda (l,k,0) \widehat{p}(\gamma) \rangle e^{2 \pi
i \lambda \alpha n} |\lambda| d\lambda \right|^2~~.
\end{eqnarray*}
Applying the Parseval formula for the interval, we thus obtain
\begin{eqnarray} \nonumber
\sum_{\gamma \in \Gamma'} | g( \gamma) |^2 & = & \alpha^{-1}
\sum_{(l,k) \in \Gamma} \int_E |\langle \widehat{f} (\lambda),
\rho_\lambda (l,k,0) \widehat{p}(\gamma) \rangle|^2 |\lambda|^2
d\lambda \\ \nonumber & = & \alpha^{-1} \int_E \sum_{(l,k) \in
\Gamma} |\langle \Phi_\lambda, T_{\lambda l} M_{k} |\lambda|^{1/2}
\mathbf{h}_\lambda \rangle e^{\pi i \lambda l k} |^2 ~|\lambda| d\lambda\\
\label{eqn:norm_samp} & = & \alpha^{-1} \int_E \sum_{(l,k) \in
\Gamma} |\langle \Phi_\lambda, T_{|\lambda|l} M_{k}
|\lambda|^{1/2} \mathbf{h}_\lambda \rangle |^2 ~|\lambda| d\lambda
\end{eqnarray} where the last equation used (\ref{eqn:scal_HS})
to express the Hilbert-Schmidt scalar products as scalar products
of vector-valued functions, as well as symmetry of $\Gamma$ to
replace $\lambda$ by $|\lambda|$.

On the other hand, by the Plancherel formula, we find
\begin{eqnarray*}
 \| g \|_2^2 & = & \int_E \| \widehat{f}(\lambda) \|_{HS}^2
 |\lambda| d\lambda \\ & = &
 \int_E \| \Phi_\lambda \|_{{\rm
 L}^2(\mathbb{R},\mathbb{C}^{d(\lambda)})}^2 |\lambda| d\lambda ~~.
\end{eqnarray*}

Hence the lower sampling estimate yields
\begin{equation} \label{eqn:ineq_I}
A \int_E \| \Phi_\lambda \|_{{\rm
 L}^2(\mathbb{R},\mathbb{C}^{d(\lambda)})}^2 |\lambda| d\lambda \le
\alpha^{-1} \int_E \sum_{(l,k) \in \Gamma} |\langle \Phi_\lambda,
T_{|\lambda|l} M_{k} |\lambda|^{1/2} \mathbf{h}_\lambda \rangle
|^2 ~|\lambda| d\lambda
\end{equation}
Since this inequality holds true for all Borel sets $E$ of
diameter at most $1/\alpha$, it has to hold pointwise a.e. for the
integrands, i.e. after shifting constants:
\begin{equation} \label{eqn:frame_fibre}
 \alpha |\lambda|^{-1} A \| \Phi_\lambda \|_{{\rm
 L}^2(\mathbb{R},\mathbb{C}^{d(\lambda)})}^2 \le
\sum_{(l,k) \in \Gamma} |\langle \Phi_\lambda, T_{|\lambda|l}
M_{k}
 \mathbf{h}_\lambda \rangle |^2 ~~(\mbox{a.e.}
\lambda)~~.\end{equation} This is already quite close to the
desired lower frame estimate, except that it holds on a set of
$\lambda$'s which may depend on the choice of $f$ (or
equivalently, on the field $(\Phi_\lambda)_{\lambda \in
\mathbb{R}^*}$).

The next step is to establish (\ref{eqn:frame_fibre}) for all $f
\in PW(\mathbb{H})$ and all $\lambda$ in a set with complement of
measure zero, independent of $f$. For this purpose we pick a
sequence $(f_n)_{n \in \mathbb{N}}$ with dense span in
$PW({\mathbb{H}})$, and obtain a set $\Omega \subset [-1,1]$ with
complement of measure zero such that (\ref{eqn:frame_fibre}) holds
for all $\lambda \in \Omega$ and all $f$ in the $(\mathbb{Q} +
i\mathbb{Q})$-span of $(f_n)_{n \in \mathbb{N}}$. But then the
$(\mathbb{Q} + i\mathbb{Q})$-span of $(\widehat{f}(\lambda))_{n
\in \mathbb{N}}$ is dense in $HS({\rm L}^2(\mathbb{R})) \circ
\widehat{P}_\lambda$, for all $\lambda$ in a Borel set $\Omega'
\subset [-1,1]$ with complement of measure zero. Hence for all
$\lambda \in \Omega \cap \Omega'$, the frame estimate holds on a
dense subset of $HS({\rm L}^2(\mathbb{R})) \circ
\widehat{P}_\lambda$, which is sufficient.

Thus we have finally established (\ref{eqn:frame_fibre}) for
almost all $\lambda \in [-1,1]$, and all $\Phi_\lambda \in {\rm
 L}^2(\mathbb{R},\mathbb{C}^{d(\lambda)})$. The same argument
applies to show the upper estimate with constant $ \alpha
|\lambda|^{-1} B$.

Now, using $\mathbf{h}_\lambda = D_{|\lambda|^{1/2}}
\mathbf{h}^{d(\lambda)}$, and the relations
\[
 M_\xi D_b = D_b M_{b \xi}~~,~~T_x D_b = D_b T_{b^{-1}x}
\] we find that
\[
T_{|\lambda|l} M_{k} \mathbf{h}_\lambda =  D_{|\lambda|^{1/2}}
\left( T_{|\lambda|^{1/2}l} M_{|\lambda|^{1/2}k}
\mathbf{h}^{d(\lambda)} \right)~~.
\] Since the image of a frame under a unitary map is a frame with identical constants, we finally obtain that
$(T_{|\lambda|^{1/2}l} M_{|\lambda|^{1/2}k}
\mathbf{h}^{d(\lambda)})_{(l,k) \in \Gamma}$ is a frame, for
almost all $|\lambda| < 1$. Now part $(b)$ follows from $d \le
d(\lambda)$ for $\lambda < \frac{1}{2d+1}$.

For the converse direction observe that by assumption on $\alpha$,
all Plancherel transforms of elements of $PW({\mathbb{H}})$ are
supported in an interval of length $1/\alpha$. Hence
(\ref{eqn:norm_samp}) holds for all $g \in PW({\mathbb{H}})$,
where this time $E = [-1,1]$, and the field
$(\Phi_\lambda)_{\lambda}$ corresponds to the Plancherel transform
of $g$. But then $(b) \Rightarrow (a)$ is immediate.

The proof that the ``almost everywhere'' contained in the
statement can be omitted for lattices relies on semi-continuity
properties of the frame bounds.

 For any unit vector
$\mathbf{f} \in {\rm L}^2(\mathbb{R},\mathbb{C}^d)$, consider the
function
\[ \Theta_{\mathbf{f}}: ]0,\frac{1}{2d+1}[ \ni \lambda \mapsto \sum_{(l,k) \in |\lambda|^{1/2}  \mathcal{M} (\mathbb{Z}^2)}
|\langle \mathbf{f}, T_{l} M_{k} |\lambda|^{1/2} \mathbf{h}^d
\rangle |^2 ~~.\] We compute
\begin{eqnarray*}
 \Theta_{\mathbf{f}} (\lambda) & = & |\lambda| \sum_{(l,k) \in |\lambda|^{1/2} \mathcal{M} (\mathbb{Z}^2)} |\langle \mathbf{f}^m, T_{l}
 M_{k}
\mathbf{h}^d \rangle |^2 \\
 & = & |\lambda| \sum_{(l,k) \in |\lambda|^{1/2} \mathcal{M} (\mathbb{Z}^2)} |\sum_{i=0}^{d}
 \langle f_i, T_{ l} M_{k} h_{i} \rangle |^2 \\
 & = & |\lambda| \sum_{(l,k)} \sum_{i,j}  \langle f_i, T_{l} M_{k} h_{i}
 \rangle  \overline{ \langle f_j, T_{l} M_{k} h_{j}
 \rangle} \\
 & = & |\lambda| \sum_{i,j} \sum_{(l,k)} \langle f_i, \langle f_j, T_{l} M_{k} h_{j}
 \rangle T_{l} M_{k} h_{i} \rangle \\
 & = & |\lambda| \sum_{i,j} \langle f_i, \sum_{(l,k)} \langle f_j, T_{l} M_{k} h_{j}
 \rangle T_{l} M_{k} h_{i} \rangle \\
 & = & |\lambda| \sum_{i,j} \langle f_i, S_{h_{i},h_{j}; |\lambda|^{1/2} \mathcal{M}} f_j \rangle ~~.
\end{eqnarray*}
Here we used the linear operator $S_{g_1,g_2; \mathcal{N}}$
associated to functions $g_1,g_2$ and an invertible matrix
$\mathcal{N}$, defined by
\[ S_{g_1,g_2; \mathcal{N}} (f) = \sum_{(l,k) \in \mathcal{N} (\mathbb{Z}^2)} \langle f, T_{l} M_{k} g_1
\rangle T_{l} M_{k} g_2 ~~.\] \cite[Theorem 3.6]{KaFe} states that
$S : M^1(\mathbb{R}) \times M^1(\mathbb{R}) \times {\rm
GL}(2,\mathbb{R}) \to \mathcal{B}(L^2(\mathbb{R}))$ is continuous,
where the right-hand side denotes the space of bounded operators
endowed with the norm topology, and $M^1(\mathbb{R})$ is the
Feichtinger algebra; see e.g. \cite{Gr} for a definition and basic
properties. Now the inclusion $\mathcal{S}(\mathbb{R}) \subset
M^1(\mathbb{R})$ entails that the map $\lambda \mapsto \langle
f_i, S_{h_i,h_j; |\lambda|^{1/2} \mathcal{M}} f_j \rangle$ is
continuous, for all $0 \le i,j \le d$, and then
$\Theta_{\mathbf{f}}$ is continuous.

Next consider the map associating to each $\lambda$ the optimal
upper frame bound, given by
\[
 B_{opt} : ]0,\frac{1}{2d+1}[ \ni \lambda \mapsto \sup_{\| \mathbf{f} \|=1} \Theta_{\mathbf{f}} (\lambda) ~~.\]
The supremum is always finite: By \cite[Corollary 6.2.3]{Gr}, the
frame operator of a one-dimensional window in the Schwartz class
is always bounded. Hence the upper frame bound also exists in the
vector valued case, by (\ref{eqn:sf_fr_bd}).

 As the supremum of a family of continuous functions, $B_{opt}$ is lower
semi-continuous, and then $\lambda \mapsto B_{opt} (\lambda)
|\lambda|^{-1/2}$ is lower semi-continuous as well. We already
know that the latter map is bounded from above by $\alpha B$ on
subset of $]0,\frac{1}{2d+1}[$ with complement of measure zero.
This subset is dense, hence lower semi-continuity implies $B_{opt}
(\lambda) |\lambda|^{-1/2} \le \alpha B$ on the whole interval.

The analogous reasoning, replacing lower by upper semi-continuity,
applies to the lower frame bound, and we are done.
\end{prf}

{\bf Proof of Theorem \ref{thm:main2}.} Fix $d \in \mathbb{N}_0$.
Suppose that $\mathcal{M}$ is given with $\| \mathcal{M} \| <
C_{\mathbb{H}}/\sqrt{2d+1}$. Let $K = {\mathcal{M}}
([-0.5,0.5)^2)$, then $K$ is a complement of $\mathcal{M}
(\mathbb{Z}^2)$ in $\mathbb{R}^2$, contained in a ball of radius
$r_0 = \| \mathcal{M} \|$, and with measure $|{\rm
det}(\mathcal{M})|$. Moreover, by choosing $\alpha>0$ small
enough, the set $K' = K \times [-\alpha/2, \alpha/2)$ is contained
in a ball of radius $r_0 + \epsilon < C_{\mathbb{H}}/\sqrt{2d+1}$
(with respect to the quasi-norm on $\mathbb{H}$). In addition,
$K'$ is a complement of $\Gamma' = \mathcal{M} (\mathbb{Z}^2)
\times \alpha \mathbb{Z}$: For any $(p,q,t) \in \mathbb{H}$, there
exist unique $(p_1,q_1) \in \Gamma$ and $(p_2,q_2) \in K$ with
$(p_1+p_2,q_1+q_2) = (p,q)$, and finally unique $l \in \mathbb{Z}$
and $s \in [-\alpha/2,\alpha/2)$ with $s+\alpha l = t - (p_1
q_2-p_2 q_1)/2$. But these choices imply $(p_1,q_1,\alpha l)
(p_2,q_2,s) = (p,q,t)$.

We will apply Theorem \ref{thm:sampl} to dilated copies of
$\Gamma'$. For this purpose, let for $a>0$ and $(p,q,t) \in
\mathbb{H}$, $\delta_a(p,q,t) = (ap,aq,a^2t)$. It is easy to check
that this defines a group automorphism $\delta_a$ fulfilling
$|\delta_a(p,q,t)| = a |(p,q,t)|$. Hence $\delta_a(\Gamma')$ is a
quasi-lattice, with complement $\delta_a(K)$ contained in a ball
of radius $a(r_0 + \epsilon)$. Hence, for any $a <
C_{\mathbb{H}}/(r_0+\epsilon)$, the sampling theorem provides the
estimate \begin{eqnarray*}
 \frac{1}{a^4 |{\rm
det}(\mathcal{M})| \alpha} \left( 1- \frac{a(r_0 +
\epsilon)}{C_{\mathbb{H}}} \right)^2 \| f \|_2^2  & \le &
\sum_{\gamma \in \delta_a(\Gamma')} |f(\gamma)|^2  \\ & \le &
\frac{1}{a^4 |{\rm det}(\mathcal{M})| \alpha} \left( 1+
\frac{a(r_0 + \epsilon)}{C_{\mathbb{H}}} \right)^2 \| f \|_2^2~~.
\end{eqnarray*}
An application of Lemma \ref{lem:WH_transfer} then yields, for all
$\lambda < \frac{1}{2n+1}$, that $\mathcal{G}(\mathbf{h}^d, a
|\lambda|^{1/2} \mathcal{M}(\mathbb{Z}^2))$ is a frame with bounds
$\frac{1}{a^2 |\lambda| ~|{\rm det}(\mathcal{M})|} \left( 1 \mp
\frac{a(r_0 + \epsilon)}{C_{\mathbb{H}}} \right)^2$. Letting $a^2
|\lambda|=1$ provides a lower frame bound for
$\mathcal{G}(\mathbf{h}^d, \mathcal{M}(\mathbb{Z}^2))$ given by
\[
 \sup \left\{ \frac{1}{|{\rm det}(\mathcal{M})|} \left( 1 -
\frac{a(r_0 + \epsilon)}{C_{\mathbb{H}}} \right)^2 ~:~ \sqrt{2d+1}
\le a < C_{\mathbb{H}}/(r_0+\epsilon) \right\}~~.
\] Observe that the restriction $a \ge \sqrt{2d+1}$ is imposed by $|\lambda| \le \frac{1}{2d+1}$.
By monotonicity, the supremum is $ \frac{1}{|{\rm det}(\mathcal{M})|} \left( 1 -
\frac{\sqrt{2d+1}(r_0 + \epsilon)}{C_{\mathbb{H}}} \right)^2$.
Sending $\epsilon$ to zero provides the lower estimate of Theorem
\ref{thm:main2}. The upper estimate is obtained in the same
fashion. \hfill $\Box$

\section*{Concluding remarks}

It is a standard observation that the construction of Gabor frames
is equivalent to the discretization of the inversion formula
associated to a certain discrete-series representation of the
so-called {\em reduced Heisenberg group} $\mathbb{H}_r$, which is
the quotient of $\mathbb{H}$ by a discrete central subgroup; see
e.g. \cite{Gr}. The proof of Theorem \ref{thm:main2} shows that in
working with $\mathbb{H}$ one needs to deal with a fair amount of
additional technical details (in particular due to the occurrence
of direct integrals), that one avoids by considering
$\mathbb{H}_r$. The benefit of this approach lies in the fact that
a single sampling estimate, namely (\ref{eqn:samp_est}), gives
rise to a whole family of Gabor frame estimates, namely
(\ref{eqn:sf_fr_bd}), valid for all $d \ge 0$.

The main results of this paper provide rather intuitive asymptotic
estimates for Gabor frame bounds. A major drawback of these
estimates is that they involve unknown constants. A ``formula''
for $C_{\mathbb{H}}$ is given in \cite{FuGr}, involving operator
norms for differential operators on $PW(\mathbb{H})$ as well as a
Sobolev constant for the unit ball ${\mathbb{H}}$; the argument is
very similar to the estimate of the constant $C_{\mathbf{f}}$ in
the proof of Theorem \ref{thm:main1}. While rough estimates for
the differential operators should be obtainable from the
Plancherel transform, which decomposes the differential operators
as well as $PW(\mathbb{H})$, we are not aware of a reasonable
estimate for the Sobolev constant for $\mathbb{H}$. In any case,
we stress that the constant $C_{\mathbb{H}}$ in the sampling
theorem is the same as in Theorem \ref{thm:main2}; this was the
chief motivation for picking the matrix norm
(\ref{eqn:defn_matr_norm}).

For single Hermite functions, the results obtained here compare in
an interesting way with recent results due to Gr\"ochenig and
Lyubarskii. Using complex analysis methods, they obtained the
following statement \cite[Theorem 3.1]{GrLy}:
\begin{thm} \label{thm:grly}
 If $|\det{\mathcal{M}}|<(d+1)^{-1}$, then
 $\mathcal{G}(h_d,{\mathcal{M}})$ is a frame for ${\rm
 L}^2(\mathbb{R})$.
\end{thm}

For the isotropic case, i.e., ${\mathcal M} = a \cdot {\rm Id}$,
this result provides a criterion that is very close to our Theorem
\ref{thm:main2}: Any $a$ below a threshold $\sim n^{-1/2}$
guarantees a frame. In the general case however, Theorem
\ref{thm:grly} is much more widely applicable: At the same time
$|\det{\mathcal{M}}|$ can be made arbitrarily small and $\|
\mathcal{M} \|$ arbitrarily large.

On the other hand, Theorem \ref{thm:grly} does not provide frame
bound estimates, and it only applies to the scalar-valued case.

Let us finally comment on possible generalizations. The first
possible extension consists in replacing $\mathbb{R}$ by
$\mathbb{R}^n$, i.e. studying vector-valued Gabor frames in ${\rm
L}^2(\mathbb{R}^n;\mathbb{C}^d)$. One now considers the
$2n+1$-dimensional Heisenberg group $\mathbb{H}_n$. This is a
stratified Lie group, possessing a sub-Laplacian, Paley-Wiener
space and, finally, a sampling theorem \cite{FuGr}. As for the
one-dimensional case, the spectral decomposition of the
sub-Laplacian involves Hermite functions, and an adaptation of the
arguments for $\mathbb{H}$ should be a straightforward task,
somewhat aggravated by additional bookkeeping.

A second, more interesting but also more challenging type of
generalization concerns the sampling sets, which could also be
irregular. There already exists an irregular sampling theorem for
$\mathbb{H}$, however, in the transfer of the associated sampling
estimates to Gabor frame estimates, we are crucially relying on
the lattice structure of the sampling set. In this context, the
key result is the continuity statement \cite[Theorem 3.6]{KaFe},
and the proof of this result makes full use of Gabor theory
developed for lattices.

As a result, we can currently only prove statements of the
following form: For all $d \in \mathbb{N}_0$ and all uniformly
discrete and uniformly dense sets $\Gamma \subset \mathbb{R}^2$
there exists a range $(0,a_d)$ of dilation parameters such that
${\mathcal{G}}(\mathbf{h}^d,a \Gamma)$ is a frame, for almost all
$a \in (0,a_d)$, including an estimate of the frame bounds.
Moreover, the threshold $a_d$ is of the order $d^{-1/2}$.

\section*{Acknowledgements}
I thank Karlheinz Gr\"ochenig and Norbert Kaiblinger for
interesting discussions and a preliminary version of \cite{GrLy}.

\end{document}